# Robust Transmission Constrained Unit Commitment: A Column Merging Method

Xuan Li [1], Qiaozhu Zhai [1*], Xiaohong Guan [1,2]

[1] Systems Engineering Institute, MOE KLINNS Lab, Xi'an Jiaotong University, Xi'an, China
[2] Center for Intelligent and Networked Systems, Department of Automation, Tsinghua University, Beijing, China
*qzzhai@sei.xjtu.edu.cn

**Abstract:** With rapid integration of power sources with uncertainty, robustness must be carefully considered in the transmission constrained unit commitment (TCUC) problem. The overall computational complexity of the robust TCUC methods is closely related to the vertex number of the uncertainty set. The vertex number is further associated with 1) the period number in the scheduling horizon as well as 2) the number of nodes with uncertain injections. In this paper, a column merging method (CMM) is proposed to reduce the computation burden by merging the uncertain nodes, while still guaranteeing the robustness of the solution. By the CMM, the transmission constraints are modified, with the parameters obtained based on an analytical solution of a uniform approximation problem, so that the computational time is negligible. The CMM is applied under a greedy-algorithm based framework, where the number of merged nodes and the approximation error can be well balanced. The CMM is designed as a preprocessing tool to improve the solution efficiency for robust TCUC problems and is compatible with many solution methods (like two-stage and multi-stage robust optimization methods). Numerical tests show the method is effective.

## 1. Introduction

The transmission constrained unit commitment (TCUC) problem is a fundamental problem in power system scheduling, planning, and electricity market clearing [1]. The goal of a typical TCUC problem is to find the optimal on/off (UC) decisions and the generation level (economic dispatch, ED) decisions of the thermal units, such that the total costs are minimized while satisfying all individual unit constraints and system-wide constraints including the transmission constraints [2]. The TCUC problem is usually modeled as a deterministic mixed-integer programming problem with load or contingency uncertainties handled by power reserves [3].

The commonly adopted methods to solve the deterministic TCUC problem are Lagrangian relaxation-based methods [4,5] and the mixed integer programming-based (MIP) methods that take advantage of the commercial solver like CPLEX and GUROBI [6,7]. For large scale TCUC problems, many techniques have been proposed to improve the computational efficiency [8], for example: [9] proposes a decomposition and coordination approach based on surrogate Lagrangian relaxation method to solve the large-scale UC problems with combined cycle units; [6,10] propose computational efficient mixed integer linear programming (MILP) formulation for UC problems with one and two sets of binary variables, respectively; [11] proposes a tight locally ideal formulation for the UC problem with piecewise linear cost functions; [12,13] propose several heuristic frameworks to obtain the near-optimal UC solutions quickly; [14] proposes a sufficient analytical condition, by which most of the redundant transmission constraints can be rapidly identified without solving any optimization problem.

The integration of large-scale renewable power sources imposes great challenges to power system operation. The volatility and intermittence of the uncertain power output must be considered in the UC problem [15]. For this purpose, stochastic programming methods, like [16-18], have been adopted, which represent the uncertainty by scenarios or scenario-trees generated with the probability distribution information of power injections of renewable energy sources. To reduce the scenario number, scenario reduction methods, like [19], are usually adopted. Chance-constrained programming methods, like [20,21], also utilize the probability distribution information and guarantee the constraints are satisfied with some probability.

However, the detailed probability information is often hard to obtain [15], and the above methods do not guarantee the robustness of the UC decisions. Robustness is crucial in actual operations, which means the UC solution of the problem can immunize against all possible realizations of the uncertainty in some predefined uncertainty set. It has been pointed out in [22,23] that a given set of UC solution is robust if and only if it is feasible for the vertex scenarios of the (polyhedral) uncertainty set.

For this purpose, the two-stage robust optimization methods like [22-25] iteratively add the vertex scenarios into the master problem. If the first-stage UC decisions cannot immunize against some specific vertex scenarios, the "worst" ones can be obtained by solving the feasibility check subproblems. However, solving the feasibility check subproblem can be hard. Many methods (reviewed in [26]) have been applied to solve the subproblem, the upper bound of computational requirements is still proportional to the vertex number of uncertainty set. To solve the feasibility check problem, parallel computing can be applied [27], but the essential difficulty in solving the feasibility check subproblem is not reduced. Also, the nonanticipativity of the dispatch solution is not well respected by these methods and can be a significant issue in real operations [28,29].

Multi-stage robust optimization methods [28,30], on the other hand, guarantee both the robustness and nonanticipativity of the solutions, by using the affine policies that replace the unit dispatch decision variables with affine functions of the uncertain power injections. The "worst" scenarios can be checked and added separately for different constraints, so the computational burden for solving feasibility check subproblems is significantly reduced. However, the affine policies are only an approximation to the full adaptive policies, so the feasible region of the UC problem is reduced, and the optimality of the solution is sacrificed.

Unlike the robust methods (both two-stage and multi-stage) that add the scenarios iteratively into the master problem, ASF-MILP method [29] guarantees the robustness and the nonanticipativity of the solution by constructing scenarios based on the structural information of the uncertainty set. The scenario number is proportional to the vertex number of the uncertainty set in a single period (instead of the vertex number of the whole uncertainty set). However, the scenarios are still too many, so the method may not be able to deal with the system with a large number of nodes with uncertainty.

It is observed that all of the above methods have the scenario-based structure in some sense. And due to the consideration of the transmission constraints, the scenario number is closely related to the vertex number of the uncertainty set. The vertex number is proportional to $2^{MT}$ with $M$ is node number with uncertain power injection/load, and $T$ is the number of time periods of the scheduling horizon.

To reduce the scenario number, in this paper, a column merging method (CMM) is proposed. The idea is to merge the $M$ uncertain nodes into $K$ nodes by using formulation transformation ($K$ is an integer and $K<M$) and thus reduce the vertex number from $2^{MT}$ to $2^{KT}$. For example, if $M-K=3$ then the vertex number is reduced to $1/2^{3T}=1/8^T$ of the original one. Meanwhile, the CMM still guarantees the solution's robustness by modifying the transmission constraints based on the errors introduced by the merging operations. The parameters of the modified constraints are obtained based on an analytical solution of a uniform approximation problem, so the computational time for obtaining the modified constraints is negligible. The CMM is applied under a greedy-algorithm based framework, where the number of merged nodes and the approximation error can be well balanced.

It must be noticed that CMM aims to reduce the formulation scale of TCUC with uncertain power (load) injections, not to solve the TCUC problem itself. The CMM is designed as a preprocessing tool and is compatible with the TCUC solution method in which the vertex number of uncertainty set is the main reason for the computational burden.

The remainder of this paper is organized as follows: In section 2, the robust TCUC problem is formulated. The column merging method is presented in section 3. In section 4, the CMM is applied and tested on two test systems, together with the ASF-MILP method and a general decomposition method, respectively. Section 5 concludes the paper.

## 2. Robust TCUC

Suppose a system has $I$ thermal units, $M$ loads, $L$ transmission lines, and is scheduled within a $T$-period horizon. Without loss of generality, the nodal net loads $\boldsymbol{d} \in \mathbb{R}^{M \times T}$ are taken as the uncertainty. It is assumed that (1) holds for every load $m$ in every period $t$.

$$d_{m,t} \in [\underline{d}_{m,t}, \overline{d}_{m,t}]; \forall m,t \qquad (1)$$

In (1), the lower and upper bounds, $\underline{d}_{m,t}, \overline{d}_{m,t} \in \mathbb{R}$, can be obtained based on forecast information. The whole uncertainty set is given by (2).

$$D = \{\boldsymbol{d} \in \mathbb{R}^{M \times T} \mid \underline{\boldsymbol{d}} \leq \boldsymbol{d} \leq \overline{\boldsymbol{d}}; f(\boldsymbol{d}) \leq \gamma\} \qquad (2)$$

Where $f: \mathbb{R}^{M \times T} \to \mathbb{R}$ is an affine function, and $f(\boldsymbol{d}) \leq \gamma$ ($\gamma \in \mathbb{R}$) represents the budget constraint.

The general robust TCUC problem can be formulated as (3)-(7):

$$\min_{z,p} S(z) + F(p) \qquad (3)$$

$$\text{s.t. } \mathbf{1}^T \boldsymbol{p}_t - \mathbf{1}^T \boldsymbol{d}_t = \mathbf{0}; \forall \boldsymbol{d} \in D, \forall t \qquad (4)$$

$$-\boldsymbol{F} \leq \boldsymbol{\Gamma}^U \boldsymbol{p}_t - \boldsymbol{\Gamma}^D \boldsymbol{d}_t \leq \boldsymbol{F}; \forall \boldsymbol{d} \in D, \forall t \qquad (5)$$

$$\boldsymbol{p} \in Y(z) \qquad (6)$$

$$z \in X, z \text{ binary} \qquad (7)$$

In this problem, the decision variables are the unit commitment variables $z \in \{0,1\}^{I \times T}$ and the dispatch variables $\boldsymbol{p} \in \mathbb{R}^{I \times T}$. The objective (3) is to minimize the total operational cost including start-up costs and fuel costs. (4) is the power balance constraint, where $\boldsymbol{p}_t$ and $\boldsymbol{d}_t$ are the $t$-th column of matrices $\boldsymbol{p}$ and $\boldsymbol{d}$, respectively. (5) represents the transmission constraints, where $\boldsymbol{\Gamma}^U \in \mathbb{R}^{L \times I}$ and $\boldsymbol{\Gamma}^D \in \mathbb{R}^{L \times M}$ are matrices of power transmission distributed factors (PTDF), and $\boldsymbol{F} \in \mathbb{R}^{L \times 1}$ is the vector of transmission limits. Both (4) and (5) should be satisfied for all possible uncertainty realizations, i.e. $\forall \boldsymbol{d} \in D$. (6) represents constraints only for the dispatch variables, including generation capacity constraints, ramping constraints, etc. Parameters of these constraints are determined with given UC decision $z$. (7) represents the constraints for UC variables such as minimum on/off constraints, must on/off constraints. More detailed formulations can be found in [10].

It should be noted that in the real operation, the uncertainty is unfolded sequentially over time, and the dispatch decisions are made accordingly. These facts are closely related with the *nonanticipative constraints* (see [28,29] for more details) and they are not explicitly formulated by (3)-(7). Equations (3)-(7) are only the basic formulation to present the critical constraints of the robust TCUC problem. The primary results of this paper are solely based on the structure of constraints (4)-(5), and are valid when the nonanticipativity of the dispatch decisions is considered.

## 3. Column Merging Method

In solving the robust TCUC problems, the computational burden is closely related to the vertex number of the uncertainty set: One simple fact is that, if the UC decisions can accommodate all the vertex scenarios of the uncertainty set, then the robustness can be guaranteed [29]. Take the box uncertainty set $D$ without budget constraints as an example: the uncertainty set has total $2^{MT}$ vertices. Therefore, with the general decomposition approach [22], the upper bound of the number of iterations is up to $2^{MT}$ in extreme cases, i.e. $O(2^{MT})$ (the number of iterations is generally much less than $2^{MT}$). Although ASF-MILP method proposed in [29] successfully guarantees the solution's robustness, by solving an MILP problem with near $2^M$ scenarios, problems with large $M$ would still have the heavy computational burden.

In this section, we propose a column merging method (CMM) to address this issue. The CMM aims to merge the columns in the matrix related to the uncertainty variables ($\mathbf{\Gamma}^D$ in (5)) and modify the constraints accordingly. With merged columns, dimensions and vertex number of the uncertainty set are reduced, and the computational burden for solving the robust TCUC problems are thus reduced.

### 3.1. Main idea

For the ease of description, equations (1),(4),(5) are rewritten in detailed form as (8)-(10):

For the robust TCUC problem, at each period $t$, we have:

$$\underline{d}_{m,t} \leq d_{m,t} \leq \overline{d}_{m,t}; \forall m \tag{8}$$

$$\sum_{i=1}^{I} p_{i,t} = \sum_{m=1}^{M} d_{m,t} \tag{9}$$

$$-F_l \leq \sum_{i=1}^{I} \Gamma^U_{l,i} p_{i,t} - \sum_{m=1}^{M} \Gamma^D_{l,m} d_{m,t} \leq F_l; \forall l \tag{10}$$

The main idea of CMM is described as follows:

First, divide the set of $M$ net loads into $K$ mutually exclusive subsets ($K \leq M$) based on some scheme (to be determined) as (11),

$$\{1,\cdots,M\} = \{h_1^{(1)},\cdots,h_{n_1}^{(1)}\} \cup \cdots \cup \{h_1^{(K)},\cdots,h_{n_K}^{(K)}\} \tag{11}$$

where $n_1 + n_2 + ... + n_K = M$.

Let

$$\tilde{d}_{k,t} = \sum_{j=1}^{n_k} d_{h_j^{(k)},t}; \forall k \in \{1,...,K\} \tag{12}$$

With (12), all $d_{m,t}$ in (8)-(9) can now be replaced by $\tilde{d}_{k,t}$. While constraint (10) is further modified as follows:

Approximate $\sum_m \Gamma^D_{l,m} d_{m,t}$ in (10) by:

$$\sum_{m=1}^{M} \Gamma^D_{l,m} d_{m,t} \approx \sum_{k=1}^{K} \left( \alpha_{l,k} \tilde{d}_{k,t} + \beta_{l,k} \right) \tag{13}$$

Where $\boldsymbol{\alpha} \in \mathbb{R}^{L \times K}$ and $\boldsymbol{\beta} \in \mathbb{R}^{L \times K}$ are approximation parameters to be determined. Denote the maximum approximation error (of all possible $d$s that satisfy (8)) of subset $k$ in (11) as $\varepsilon_{l,k}$. Then we have:

$$\left| \sum_{j=1}^{n_k} \Gamma^D_{l,h_j^{(k)}} d_{h_j^{(k)},t} - \left( \alpha_{l,k} \tilde{d}_{k,t} + \beta_{l,k} \right) \right| \leq \varepsilon_{l,k} \tag{14}$$

With (11)-(14), (8)-(10) are now replaced by (15)-(17):

$$\sum_{j=1}^{n_k} \underline{d}_{h_j^{(k)},t} \leq \tilde{d}_{k,t} \leq \sum_{j=1}^{n_k} \overline{d}_{h_j^{(k)},t}; \forall k \tag{15}$$

$$\sum_{i=1}^{I} p_{i,t} = \sum_{k=1}^{K} \tilde{d}_{k,t} \tag{16}$$

$$-F_l + \sum_{k=1}^{K} \varepsilon_{l,k} \leq \sum_{i=1}^{I} \Gamma^U_{l,i} p_{i,t} \\ - \sum_{k=1}^{K} \left( \alpha_{l,k} \tilde{d}_{k,t} + \beta_{l,k} \right) \leq F_l - \sum_{k=1}^{K} \varepsilon_{l,k}; \forall l \tag{17}$$

It can be seen from (15)-(17) that, the $M$ nodes with uncertainty are merged into $K$ nodes, and the related constraints are modified accordingly. The relationship between the original constraints and modified constraints can be summarized in theorem 1:

**Theorem 1**: Suppose $d$ lies in the box region defined by (8), $\tilde{d}$ is determined based on (12). If ED decision $p$ satisfies (16)-(17), then it must also satisfy (9)-(10).

The result can be directly obtained from (14), and the proof is omitted in this paper. Theorem 1 essentially means (15)-(17) are sufficient conditions of (8)-(10). In other words, if a set of UC solution $z$ is a feasible solution of problem (3)(6)(7)(16)(17) with any $\tilde{d}$ satisfying (15), then it must also be a feasible solution of problem (3)(6)(7)(9)(10) with any $d$ satisfying (8). Therefore, with merged nodes and modified constraints, the solution's robustness is not sacrificed.

Now, two questions are remained to be solved: 1) How to determine the approximation parameters $(\alpha, \beta)$ and error $\varepsilon$ in (17)? 2) How to divide the nodes into groups, as in (11)? The answers are presented in the following two subsections.

### 3.2. Calculation of Approximation Parameters

Approximation parameters $(\alpha, \beta, \varepsilon)$ in (17) are determined as follows:

For each $k$ ($1 \leq k \leq K$) and transmission line $l$ ($1 \leq l \leq L$), solve the following minimax optimization problem:

$$\varepsilon_{l,k} = \min_{\alpha_{l,k},\beta_{l,k}} \max_{d_{h_j^{(k)},t}} \left| \sum_{j=1}^{n_k} \Gamma^D_{l,h_j^{(k)}} d_{h_j^{(k)},t} - \alpha_{l,k} \sum_{j=1}^{n_k} d_{h_j^{(k)},t} - \beta_{l,k} \right| \tag{18}$$

s.t. $\underline{d}_{h_j^{(k)},t} \leq d_{h_j^{(k)},t} \leq \overline{d}_{h_j^{(k)},t}$ ; $j = 1,2,\cdots,n_k$ (19)

Suppose the optimal solution is $(\alpha^*_{l,k}, \beta^*_{l,k})$ with the optimal objective value of $\varepsilon^*_{l,k}$, then let $(\alpha, \beta, \varepsilon) = (\alpha^*, \beta^*, \varepsilon^*)$.

The above procedure is, in fact, a kind of best uniform approximation [31]. To solve the problem (18)-(19), an analytical solution is given in this subsection, so the computational time is negligible:

For clarity, (18)-(19) can be equivalently formulated as the following mathematical programming problem (index $k$ and $l$ are omitted):

$$\varepsilon_0 = \min_{\alpha_0, \beta_0} \max_{x_j} \left| \sum_{j=1}^{n} \alpha_j x_j - \alpha_0 (\sum_{j=1}^{n} x_j) - \beta_0 \right| \quad (20)$$

$$\text{s.t. } 0 \leq x_j \leq \bar{x}_j ; j = 1, 2, \cdots, n \quad (21)$$

Notice that the lower bounds of variable $x$ are 0 (can be obtained by simple linear transformation). Without loss of generality, it is always assumed that:

$$\alpha_1 \leq \alpha_2 \leq \ldots \leq \alpha_n \quad (22)$$

For otherwise, the subscripts can be reset such that (22) is satisfied.

Now, a sequence is defined as follows:

$$\lambda_j = \sum_{i=1}^{j-1} \bar{x}_i - \sum_{i=j}^{n} \bar{x}_i ; j = 1, \cdots, n+1 \quad (23)$$

And there must exist a $j^*$ ($1 \leq j^* \leq n$) such that

$$\lambda_{j^*} \leq 0, \lambda_{j^*+1} \geq 0 \quad (24)$$

With the above notations, we have theorem 2:

**Theorem 2**: $(\alpha_0^*, \beta_0^*)$ obtained by (25)-(26) is the globally optimal solution to the problem (20)-(21).

$$\alpha_0^* = \alpha_{j^*} \quad (25)$$

$$\beta_0^* = \frac{1}{2} \sum_{j=1}^{n} (\alpha_j - \alpha_0^*) \bar{x}_j \quad (26)$$

Proof of the theorem 2 is presented in the Appendix.

Theorem 2 essentially gives an analytical solution to the problem (20)-(21). So, the solution time for this problem is negligible. The importance of this feature will be seen in the next subsection.

### 3.3. Algorithmic Framework

The only remaining problem now is how to determine the partition of the original nodes. The requirement for the partition should be in two aspects: 1) the number of merged nodes ($K$ in (11)). The smaller the $K$ is, the smaller scale the modified problem will have. 2) Approximation error ($\varepsilon$ in (17)): The smaller the $\varepsilon$ is, the more reduced feasible region (with more tightened transmission constraints) the modified problem will have.

Suppose the number of merged node $K$ should be no more than $K^{req}$, and approximation error should be no more than $E^{req}$. In this subsection, a framework based on greedy algorithm [32] is proposed to quickly obtain the satisfactory partition.

**(Algorithm 1)**

**Step 1**: Initialize the set of nodes as $\Lambda = \{\Lambda_1, \cdots, \Lambda_M\}$ where $\Lambda_m = \{m\}$, and the number of merged nodes as $K^{temp} = M$. Set all elements in $\varepsilon^{temp} \in \mathbb{R}^{L \times K^{temp}}$ as zeros.

**Step 2**: If $K^{temp} \leq K^{req}$, then go to step 7. Otherwise, go to step 3.

**Step 3**: Obtain $\varepsilon_l^{j_1, j_2}$ of all $j_1, j_2 \in \Lambda$ ($j_1 \neq j_2$) based on theorem 2 ((46) in the proof), where $\varepsilon_l^{j_1, j_2}$ is the approximation error of line $l$ merging all the nodes in the set $\Lambda_{j_1} \cup \Lambda_{j_2}$ into one node.

**Step 4**: find:

$$(j_1^*, j_2^*) = \underset{\Lambda_{j_1}, \Lambda_{j_2} \in \Lambda (j_1 \neq j_2)}{\arg\min} \max_l \varepsilon_l^{j_1, j_2} \quad (27)$$

**Step 5**: Add the column $\varepsilon^{j_1^*, j_2^*} \in \mathbb{R}^{L \times 1}$ to the right of $\varepsilon^{temp}$; Remove the $j_1^*$-th and $j_2^*$-th columns from $\varepsilon^{temp}$.

Let $\varepsilon^{check} = \max_l \sum_{k=1}^{K^{temp}} \varepsilon_{l,k}^{temp}$.

If $\varepsilon^{check} \geq E^{req}$, then go to step 7. Otherwise, go to step 6.

**Step 6**: Add $\Lambda_{j_1^*} \cup \Lambda_{j_2^*}$ into $\Lambda$. Remove $\Lambda_{j_1^*}$ and $\Lambda_{j_2^*}$ from $\Lambda$. $K^{temp} = K^{temp} - 1$. Go to step 2.

**Step 7:** End the algorithm with the output $\Lambda$.

**Remark 1**. In step 3, $\varepsilon_l^{j_1, j_2}$ of all $K^{temp}(K^{temp}-1)/2$ combinations of $j_1$ and $j_2$ must be obtained. Although it seems burdensome, in each iteration besides the first iteration, only $K^{temp} - 1$ combinations are calculated, and the others are those whose errors have been calculated previously and stored. In the first iteration, $M(M-1)/2$ combinations are calculated.

**Remark 2**. As analyzed in remark 1, $\varepsilon_l^{j_1, j_2}$ must be repeatedly calculated during the whole procedure. With the analytical solution of the problem (20)-(21) given in theorem 1, the calculation time of $\varepsilon_l^{j_1, j_2}$ can be negligible.

**Remark 3.** Although not stated in the algorithm, it is better to eliminate the redundant transmission constraints previously before running algorithm 1: It can be seen from step 3, with a reduced set of transmission lines, fewer times of $\varepsilon_l^{j_1, j_2}$ are calculated, and the maximum approximation error can be smaller. Moreover, the scale of the robust TCUC problem can be further reduced if the redundant transmission constraints are removed from the problem. A fast identification method proposed in [14] is adopted in this paper, by which most of the redundant transmission constraints can be identified quickly without solving any optimization problem.

**Remark 4.** With the proposed CMM method, the nodal price can still be determined in a traditional way: Suppose the robust UC solution obtained with CMM is $z'$. Then with given $z'$, the TCED problem with the original transmission constraints (without the CMM) is solved. The nodal price can be determined based on the optimal dual solutions of the TCED problem.

### 4. Numerical Results

As interpreted in section 1, CMM aims to reduce the formulation scale of TCUC with uncertain power (load) injections. CMM is not an algorithm for solving TCUC itself but a tool for accelerating any procedure for solving TCUC. It can *work together* with the solution methods of TCUC, in which the vertex number of uncertainty set is the main reason for the computational burden.

**Table 1** Column Merging Results

| K | Node Indices | Max. $\delta$ | Avg. $\delta$ |
|---|---|---|---|
| 8 | {1},{2},{3},{4},{5},{6},{7},{8} | 0.00% | 0.00% |
| 7 | {1},{2},{3},{4},**{5,7}**,{6},{8} | 9.42% | 0.97% |
| 6 | {1},{2},{3},{4},**{5,7,8}**,{6} | 11.70% | 1.41% |
| 5 | {1},**{2,3}**,{4},{5,7,8},{6} | 12.80% | 2.93% |
| 4 | **{1,2,3}**,{4},{5,7,8},{6} | 13.68% | 3.49% |
| 3 | {1,2,3},**{4,6}**,{5,7,8} | 22.83% | 5.12% |
| 2 | **{1,2,3,4,6}**,{5,7,8} | 24.08% | 5.73% |
| 1 | **{1,2,3,4,5,6,7,8}** | 41.32% | 6.83% |

Therefore, in this section, CMM is combined with two specific solution methods to see the computational improvements with different numbers of merged nodes: 1) the all-scenario-feasible MILP method [29] on a modified IEEE 118-bus system with wind power sources, and 2) a general robust decomposition method (similar to [33]) on the IEEE 118-bus systems with load uncertainty.

All numerical tests are performed with MATLAB R2015b, YALMIP toolbox [34] and GUROBI 7.5 package with default Gap settings of 0.01%, on an Intel(R) Core(TM) i7-3770 CPU @ 3.40GHz PC with 8GB RAM.

*4.1. IEEE 118-bus System with Wind Uncertainty*

The system is scheduled in a 24-hour horizon with total 54 schedulable thermal units, 179 transmission lines and 91 loads. Total 8 wind farms are added into the system on bus 14, 28, 42, 56, 70, 84, 98, 112, respectively. Each wind farm has the same 140 MW capacity, so 1120 MW in total (18.7% of the peak system load, 35.0% of the valley load). The transmission limits of all lines are expanded by 140MW. The upper bounds for all possible wind power outputs are the installed capacities, and the lower bounds are zeros.

We first remove the redundant transmission constraints by using the technique presented in [14]. Total 5368 (62.8%) of all transmission constraints (2 × 179 × 24=8544) are identified redundant. Moreover, no active transmission constraints are remaining for 42 lines. Transmission constraints of these lines are removed from the following CMM procedure.

The CMM is now applied to this system. The merging procedure is presented in Table 1, where the first column indicates the total number of merged wind nodes (denoted as *K*), in the second column, each brace represents a node, the numbers in which are the original indices of the wind farms. In the third and fourth columns, "Max. $\delta$" means the maximum relative error (to the transmission capacities) of all transmission lines and "Avg. $\delta$" is the average relative error.

$$\text{Max. } \delta = \max_l \frac{\sum_k \varepsilon_{l,k}}{F_l}, \quad \text{Avg. } \delta = \frac{1}{L}\sum_l \frac{\sum_k \varepsilon_{l,k}}{F_l}$$

It is seen from the above results that 1) the merging result for *K* is based on the result for *K*+1, this is the key feature of greedy algorithm shown in Algorithm 1. 2) As *K* decreases, the errors monotonically increase, since every merging operation would bring in new errors to the existing errors on the transmission lines.

With the results obtained by the CMM, we now test the all-scenario-feasible (ASF) MILP method [29] on this system. The ASF-MILP method has a scenario-based MILP form, with $2^K$ selected vertex scenarios, 1 expected scenario, and a set of nonanticipative constraints. Each scenario is a $K \times T$ matrix, i.e. $w^s \in \mathbb{R}^{K \times T}$. The robustness and nonanticipativity of the solution can be guaranteed.

It is noticed that the ASF-MILP method, though has the similar scenario-based form, is completely different from the traditional stochastic scenario-based method (SO): 1) the scenarios in SO are generated by sampling using probabilistic distribution, while ASF-MILP constructs scenarios based on the information of the uncertainty set. 2) ASF-MILP is able to guarantee the solution's robustness, while SO cannot. The readers are referred to [29] for more details.

In this test, the objective is to minimize the cost of the expected scenario. The time limit for the solver is set as 7200s (2 hours). The results are presented in Table 2, where "*K*" is the number of merged nodes with uncertainty as Table 1, "Time" is the total computational time of the ASF-MILP method spent by the MIP solver. "Scenarios" represents the total number of the $2^K$ selected vertex scenarios and the expected scenario in the UC problem (e.g. for *K*=8, the scenario number is $2^K + 1 = 2^8 + 1 = 257$). "Obj. Value" is the optimal objective value of the ASF-MILP.

Then we conduct the Monte Carlo simulation on the obtained UC solution to examine its economic performance. We totally generate 1000 samples with various probability distributions (each sample contains the an $8 \times T$ matrix of wind power outputs).

Each element of the $8 \times T$ matrix of wind power outputs is a random variable. In our simulation, the probability distribution of each element is different. For the ease of testing, the joint distribution function of these random variables is set as the product of all marginal distribution functions, i.e., these random variables are assumed to be independent and therefore are not correlated. It is also clarified that the proposed CMM method is still valid when the spatial/temporal correlation of wind power outputs are considered.

For each sample, a modified economic dispatch (ED) problem (considering the nonanticipativity of the ED decisions) is solved and cost of this sample is obtained. The results are also presented in Table 2, where "Avg. Cost" is the average cost of the 1000 samples.

**Table 2** Results of ASF-MILP method using CMM

| K | ASF-MILP | | | M.C. Sim. | |
|---|---|---|---|---|---|
| | Time (s) | Scenarios | Obj. Value ($\times 10^6$ \$) | Samples | Avg. Cost ($\times 10^6$ \$) |
| 8 | Time out | 257 | -- | | -- |
| 7 | 7122.7 | 129 | 1.3829 | | 1.3834 |
| 6 | 1467.3 | 65 | 1.3828 | | 1.3832 |
| 5 | 741.2 | 33 | 1.3846 | 1000 | 1.3840 |
| 4 | 127.6 | 17 | 1.3849 | | 1.3847 |
| 3 | 53.9 | 9 | 1.3853 | | 1.3852 |
| 2 | 29.1 | 5 | 1.3863 | | 1.3854 |
| 1 | 12.7 | 3 | 1.3873 | | 1.3859 |

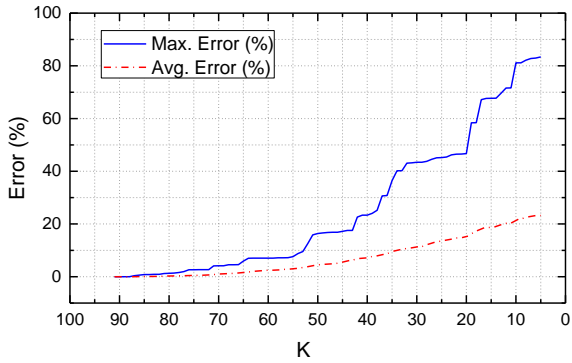

Fig. 1 Errors of CMM

It is seen from the above results that: 1) as $K$ decreases, the computational time is reduced significantly. Since the scale of the MILP problem is proportional to its scenario number, and the small-scale problem is much easier to solve than the large ones. When $K$=8 (the original problem without CMM), the problem is too large (257 scenarios) to solve within the 7200s time limit. 2) As $K$ decreases, the optimal objective value gets larger, where the optimal value of $K$=1 is 0.32% larger than the one of $K$=7. As suggested in Table 1, the column merging operation brings in errors that can reduce the feasible region of the problem, and the optimality of the solution is thus affected. However, the differences between the simulated costs are relatively small (the cost of $K$=1 is only 0.18% larger than $K$=7). It has been pointed out in theorem 1 that the robustness of the solution is not affected using the CMM, compared with such large computational improvement (the time of $K$=1 is only 0.18% of $K$=7), the sacrifice of the solution's optimality seems acceptable.

### 4.2. IEEE 118-bus System with Load Uncertainty

The basic settings of the system are the same with the one in Section 4.1, except that: 1) The 91 loads are taken as the uncertainty with 110% and 90% of the forecasted loads as the upper and lower bounds of the uncertainty set. 2) No wind farms are added into the system. 3) The transmission capacities are not expanded.

The redundant constraints are also removed before applying the CMM method. Moreover, the accumulated maximum and average errors of all transmission lines over the number of merged nodes are presented in Fig. 1.

It is seen from the above results that: 1) as $K$ decreases, the accumulated errors increase. 2) Compared to the preceding case with wind uncertainty, the errors increase relatively slow: The maximum error exceeds 10% when the 91 original loads are merged into 52 nodes, and 20% when merged into 42 nodes. From the optimization problem (20)-(21), we can see that the merging error is related to the difference between the upper and lower bounds of the uncertainty interval (i.e. $\bar{d}_{m,t} - \underline{d}_{m,t}$), and the smaller the interval is, usually the smaller the merging error will be. 3) When the 91 loads are merged into 4 nodes, the accumulated errors of at least one transmission lines exceed 100%, and the corresponding transmission constraints can never be satisfied.

The column merging results are then tested with a general robust decomposition method (similar to [33]). The method has a master-sub-problem framework: the master problem (MP) is a scenario-based MILP problem where the scenarios are added iteratively. The cost of the expected scenario is minimized in the MP. The subproblem (SP) serves as a feasibility check problem to see if the solution of the master problem is feasible for all possible scenarios within the uncertainty set. If the feasibility check is passed, then the robust solution is found. Otherwise, the optimal solution (load scenario) of the SP is added into the master problem. Then the MP is resolved, and so on.

Using CMM in the decomposition method has two benefits: 1) based on the analysis at the beginning of section 3 that CMM can reduce the vertex number of the uncertainty set and thus reduce the maximum possible number of iterations: $O(2^{MT}) \rightarrow O(2^{KT})$ ($K \leq M$). 2) It has been pointed out in [22] that if the solution of the MP is not robust, then the optimal solution of the feasibility check SP must be obtained on one of the vertices of the uncertainty set. In this test, we use extreme point based solution method (like in [23]) (where vertices are explicitly formulated by using binary variables) to solve the SP. The binary variables in the SP (after dualization and linearization, now is a MILP problem) are reduced from $MT$ to $KT$ ($M$ to $K$, if time-decoupled solution method is used), so the SP can be easier to solve after using CMM.

The results regarding the computational performance are presented in Table 3. Where the notations in the first three columns are the same with Table 2. "Iter." In the fourth column means the number of iterations used to find the robust solution, and "Time/Iter." is the average time used per iteration.

**Table 3** Computational Performance of Robust Decomposition Method using CMM

| K | Max. $\delta$ (%) | Time (s) | Iter. | Time/Iter.(s) |
|---|---|---|---|---|
| 91 | 0 | 5182.4 (100%) | 7 | 740.3 |
| 80 | 1.28% | 4359.9 (84.1%) | 11 | 396.4 |
| 70 | 4.15% | 2492.6 (48.1%) | 8 | 311.6 |
| 60 | 7.04% | 238.6 (4.6%) | 2 | 119.3 |
| 50 | 16.40% | 944.4 (18.2%) | 6 | 157.4 |
| 40 | 23.35% | 854.0 (16.5%) | 4 | 213.5 |
| 30 | 43.40% | 1353.6 (26.1%) | 7 | 193.4 |
| 20 | 46.72% | 575.0 (11.1%) | 12 | 47.9 |
| 10 | 81.15% | Infeasible | -- | -- |

**Table 4** Economy Performance of Robust Decomposition Method using CMM

| K | Max. $\delta$ (%) | Obj. V. ($\times 10^6$ $) | Sim. Cost ($\times 10^6$ $) |
|---|---|---|---|
| 91 | 0 | 1.6475 (+0.0%) | 1.6466 (+0.0%) |
| 80 | 1.28% | 1.6477 (+0.01%) | 1.6469 (+0.02%) |
| 70 | 4.15% | 1.6575 (+0.61%) | 1.6515 (+0.30%) |
| 60 | 7.04% | 1.6611 (+0.83%) | 1.6532 (+0.40%) |
| 50 | 16.40% | 1.6706 (+1.40%) | 1.6577 (+0.67%) |
| 40 | 23.35% | 1.6751 (+1.68%) | 1.6596 (+0.79%) |
| 30 | 43.40% | 1.7191 (+4.35%) | 1.6916 (+2.73%) |
| 20 | 46.72% | 1.7126 (+3.95%) | 1.6825 (+2.18%) |
| 10 | 81.15% | -- | -- |

It is seen from the results that: 1) All cases with the merged nodes (K<91) spend less time than the one with the original nodes (K=91). 2) The total time does not monotonically decrease as K decreases (with 238.6s when K=60, but 1353.6s when K=30). This is because, in the semi-enumerated decomposition approach, the number of iterations used to find the robust solution is in fact hard to estimate as the fourth column suggests (only 2 iterations when K=60). 3) No feasible solution can be found for K=10 with 81.15% accumulated maximum error.

Monte Carlo simulations are also conducted on the solutions of the decomposition method like the preceding subsection. The objective values and the average simulated costs are presented in Table 4.

It is seen from the results that 1) As more nodes are merged (K decreases), the objective value and the simulated costs increase in general. This is because the feasible region is reduced due to the introduced errors of the transmission constraints. 2) the differences between the simulated costs are smaller than differences between the objective values. 3) When 40 nodes are merged into 30, the "Max. $\delta$" turns from 23.35% to 43.40%, and both objective value and the simulated cost increase significantly (2.63% and 1.93% respectively) 3) When 30 nodes are merged into 20 (both of these cases have "Max. $\delta$" larger than 40%), the two values decrease.

## 5. Conclusions

Among the methods for the robust TCUC problem, many have the scenario-based structure. The scenario number, as well as the problem scale and overall solution difficulties, are closely related to the vertex number of the uncertainty set (the more, the harder). The CMM proposed in this paper can merge the nodes with uncertainty so that the difficulties solving the robust TCUC problems can be reduced. Meanwhile, the robustness of the solution is still guaranteed. A trade-off exists between the final number of the merged nodes and the error introduced by the merging operation. However, if the stopping criterion of the proposed greedy algorithms is chosen appropriately, the overall computational efficiency can be improved significantly while the sacrifice of the solution's optimality is negligible. The CMM is compatible with many kinds of solution methods including the ASF-MILP method and the general decomposition method tested in this paper.

**Appendix: Proof of Theorem 2**

The following auxiliary functions are useful in this part:

$$\Phi(x,\alpha_0) = \sum_{j=1}^{n} \alpha_j x_j - \alpha_0 (\sum_{j=1}^{n} x_j) \quad (28)$$

$$\Phi^{\max}(\alpha_0) = \max_{x} \Phi(x,\alpha_0) = \max_{0 \le x_j \le \bar{x}_j} \left[ \sum_{j=1}^{n} \alpha_j x_j - \alpha_0 (\sum_{j=1}^{n} x_j) \right] \quad (29)$$

$$\Phi^{\min}(\alpha_0) = \min_{x} \Phi(x,\alpha_0) = \min_{0 \le x_j \le \bar{x}_j} \left[ \sum_{j=1}^{n} \alpha_j x_j - \alpha_0 (\sum_{j=1}^{n} x_j) \right] \quad (30)$$

Theorem 2 is proved by proving the following three conclusions:

(i) $\Phi^{\max}(\alpha_0) - \Phi^{\min}(\alpha_0)$ is piecewise linear and convex for $\alpha_0 \in (-\infty, +\infty)$.

(ii) $\max_{0 \le x_j \le \bar{x}_j} \left| \Phi(x,\alpha_0) - \tilde{\beta}_0 \right| \ge \frac{1}{2} \left[ \Phi^{\max}(\alpha_0) - \Phi^{\min}(\alpha_0) \right] \quad (31)$

And " = " holds if and only if

$$\tilde{\beta}_0 = \frac{1}{2} \left[ \Phi^{\max}(\alpha_0) + \Phi^{\min}(\alpha_0) \right] \quad (32)$$

(iii) $(\alpha_0^*, \beta_0^*)$ obtained by (25)(26) is the global optimal solution to problem (20)(21).

*Proof*: (i) equation (29)-(30) can be rewritten as (33).

$$\begin{cases} \Phi^{\max}(\alpha_0) = \max_{0 \le x_j \le \bar{x}_j} \sum_{j=1}^{n} (\alpha_j - \alpha_0) x_j \\ \Phi^{\min}(\alpha_0) = \min_{0 \le x_j \le \bar{x}_j} \sum_{j=1}^{n} (\alpha_j - \alpha_0) x_j \end{cases} \quad (33)$$

Combine (22) and (33) we have the following results: If $-\infty < \alpha_0 \le \alpha_1$, then $\alpha_j - \alpha_0 \ge 0$ holds for $j = 1, 2, \cdots, n$. And thus

$$\begin{cases} \Phi^{\max}(\alpha_0) = \max_{0 \le x_j \le \bar{x}_j} \sum_{j=1}^{n} (\alpha_j - \alpha_0) x_j = \sum_{j=1}^{n} (\alpha_j - \alpha_0) \bar{x}_j \\ \Phi^{\min}(\alpha_0) = \min_{0 \le x_j \le \bar{x}_j} \sum_{j=1}^{n} (\alpha_j - \alpha_0) x_j = 0 \end{cases} \quad (34)$$

$$\Phi^{\max}(\alpha_0) - \Phi^{\min}(\alpha_0) = \sum_{j=1}^{n} \alpha_j \bar{x}_j - \alpha_0 \sum_{j=1}^{n} \bar{x}_j \quad (35)$$

If $a_{\tilde{j}} < a_0 \le a_{\tilde{j}+1}$ holds for some $\tilde{j} (1 < \tilde{j} \le n-1)$ then $a_j - a_0 \ge 0$ holds for $j = \tilde{j}+1, \tilde{j}+2, \cdots, n$ and thus

$$\begin{cases} \Phi^{\max}(\alpha_0) = \max_{0 \le x_j \le \bar{x}_j} \sum_{j=1}^{n} (\alpha_j - \alpha_0) x_j = \sum_{j=\tilde{j}+1}^{n} (\alpha_j - \alpha_0) \bar{x}_j \\ \Phi^{\min}(\alpha_0) = \min_{0 \le x_j \le \bar{x}_j} \sum_{j=1}^{n} (\alpha_j - \alpha_0) x_j = \sum_{j=1}^{\tilde{j}} (\alpha_j - \alpha_0) \bar{x}_j \end{cases} \quad (36)$$

$$\Phi^{\max}(\alpha_0) - \Phi^{\min}(\alpha_0) = \sum_{j=\tilde{j}+1}^{n} \alpha_j \bar{x}_j - \sum_{j=1}^{\tilde{j}} \alpha_j \bar{x}_j + \alpha_0 (\sum_{j=1}^{\tilde{j}} \bar{x}_j - \sum_{j=\tilde{j}+1}^{n} \bar{x}_j). \quad (37)$$

If $a_n < a_0 < +\infty$, then $a_j - a_0 \le 0$ holds for $j = 1, 2, \cdots, n$. and thus

$$\begin{cases} \Phi^{\max}(\alpha_0) = \max_{0 \le x_j \le \bar{x}_j} \sum_{j=1}^{n} (\alpha_j - \alpha_0) x_j = 0 \\ \Phi^{\min}(\alpha_0) = \min_{0 \le x_j \le \bar{x}_j} \sum_{j=1}^{n} (\alpha_j - \alpha_0) x_j = \sum_{j=1}^{n} (\alpha_j - \alpha_0) \bar{x}_j \end{cases} \quad (38)$$

$$\Phi^{\max}(\alpha_0) - \Phi^{\min}(\alpha_0) = -\sum_{j=1}^{n} \alpha_j \bar{x}_j + \alpha_0 \sum_{j=1}^{n} \bar{x}_j \quad (39)$$

Based on the above analysis it is seen that $\Phi^{\max}(\alpha_0) - \Phi^{\min}(\alpha_0)$ is piecewise linear, and the derivative of it is monotonically increasing. Conclusion (i) is hence proved.

(ii) It is seen that $\Phi^{\max}(\alpha_0)$ and $\Phi^{\min}(\alpha_0)$ are two extremal values of $\Phi(x,\alpha_0)$. Therefore, we have:

$$\max_{0 \leq x_j \leq \bar{x}_j} \left| \Phi(x, \alpha_0) - \tilde{\beta}_0 \right|$$
$$= \max \left\{ \left| \Phi^{\max}(\alpha_0) - \tilde{\beta}_0 \right|, \left| \Phi^{\min}(\alpha_0) - \tilde{\beta}_0 \right| \right\} \quad (40)$$

If $\tilde{\beta}_0 > \frac{1}{2} \left[ \Phi^{\max}(\alpha_0) + \Phi^{\min}(\alpha_0) \right]$, then:

$$\left| \Phi^{\min}(\alpha_0) - \tilde{\beta}_0 \right| = \left| \tilde{\beta}_0 - \Phi^{\min}(\alpha_0) \right|$$
$$= \left| \tilde{\beta}_0 - \frac{1}{2} \left[ \Phi^{\max}(\alpha_0) + \Phi^{\min}(\alpha_0) \right] \right.$$
$$\left. + \frac{1}{2} \left[ \Phi^{\max}(\alpha_0) - \Phi^{\min}(\alpha_0) \right] \right|$$
$$> \frac{1}{2} \left[ \Phi^{\max}(\alpha_0) - \Phi^{\min}(\alpha_0) \right] \quad (41)$$

If $\tilde{\beta}_0 < \frac{1}{2} \left[ \Phi^{\max}(\alpha_0) + \Phi^{\min}(\alpha_0) \right]$ then similarly we have:

$$\left| \Phi^{\max}(\alpha_0) - \tilde{\beta}_0 \right| > \frac{1}{2} \left[ \Phi^{\max}(\alpha_0) - \Phi^{\min}(\alpha_0) \right]. \quad (42)$$

If $\tilde{\beta}_0 = \frac{1}{2} \left[ \Phi^{\max}(\alpha_0) + \Phi^{\min}(\alpha_0) \right]$ holds then it is clear that

$$\max_{0 \leq x_j \leq \bar{x}_j} \left| \Phi(x, \alpha_0) - \tilde{\beta}_0 \right|$$
$$= \max \left\{ \left| \Phi^{\max}(\alpha_0) - \tilde{\beta}_0 \right|, \left| \Phi^{\min}(\alpha_0) - \tilde{\beta}_0 \right| \right\}$$
$$= \max \left\{ \frac{1}{2} \left[ \Phi^{\max}(\alpha_0) - \Phi^{\min}(\alpha_0) \right], \right.$$
$$\left. , \frac{1}{2} \left[ \Phi^{\max}(\alpha_0) - \Phi^{\min}(\alpha_0) \right] \right\} \quad (43)$$

Conclusion (ii) is then proved.

(iii) Combining (20)-(21), (28), and Conclusion (ii), it is clear that

$$\varepsilon_0 = \min_{\alpha_0, \tilde{\beta}_0} \max_{x_j} \left| \Phi(x, \alpha_0) - \tilde{\beta}_0 \right|$$
$$= \min_{\alpha_0, \tilde{\beta}_0} \frac{1}{2} \left[ \Phi^{\max}(\alpha_0) - \Phi^{\min}(\alpha_0) \right] \quad (44)$$

And the optimal value of $\tilde{\beta}_0$ can be obtained based on (32) if only the optimal value of $\alpha_0$ is obtained. For the $j^*$ defined by (23)-(24), it is seen from (35)(37) and (39) that

$$\frac{d}{d\alpha_0} \left( \Phi^{\max}(\alpha_0) - \Phi^{\min}(\alpha_0) \right)$$
$$= \begin{cases} \sum_{j=1}^{j^*-1} \bar{x}_j - \sum_{j=j^*}^{n} \bar{x}_j, & \text{if } \alpha_{j^*-1} < \alpha_0 < \alpha_{j^*} \\ \sum_{j=1}^{j^*} \bar{x}_j - \sum_{j=j^*+1}^{n} \bar{x}_j, & \text{if } \alpha_{j^*} < \alpha_0 < \alpha_{j^*+1} \end{cases}$$
$$= \begin{cases} \lambda_{j^*}, & \text{if } \alpha_{j^*-1} < \alpha_0 < \alpha_{j^*} \\ \lambda_{j^*+1}, & \text{if } \alpha_{j^*} < \alpha_0 < \alpha_{j^*+1} \end{cases} \quad (45)$$

Based on (24) it is known that $\Phi^{\max}(\alpha_0) - \Phi^{\min}(\alpha_0)$ is monotonically decreasing when $\alpha_{j^*-1} < \alpha_0 < \alpha_{j^*}$, and monotone increasing when $\alpha_{j^*} < \alpha_0 < \alpha_{j^*+1}$. Therefore $a_{j^*}$ is a local minimum point of $\Phi^{\max}(\alpha_0) - \Phi^{\min}(\alpha_0)$. However, the function is convex according to conclusion (i) and thus $\alpha_{j^*}$ is, in fact, the global minimum point. Together with (44) it is seen that

$$\varepsilon_0^* = \min_{\alpha_0, \tilde{\beta}_0} \frac{1}{2} \left[ \Phi^{\max}(\alpha_0) - \Phi^{\min}(\alpha_0) \right]$$
$$= \frac{1}{2} \left[ \Phi^{\max}(\alpha_{j^*}) - \Phi^{\min}(\alpha_{j^*}) \right] \quad (46)$$

and (25) is proved. (26) can be obtained by directly substituting (25) into (32).

Q.E.D